\documentclass[10pt,a4paper]{article}

\usepackage[english]{babel}
\usepackage[T1]{fontenc}
\usepackage[utf8]{inputenc}
\usepackage{arev}

\usepackage[margin=2.5cm]{geometry}
\usepackage{amsmath}
\usepackage{amssymb}
\usepackage{amsthm}

\usepackage{graphicx}
\usepackage{xcolor}
\usepackage{adjustbox}
\usepackage{tikz}
\usepackage{pgfplots}
\pgfplotsset{compat=1.15}
\usetikzlibrary{calc,patterns,angles,quotes}

\newtheorem{theorem}{Theorem}[section]

\theoremstyle{definition}
\newtheorem{definition}{Definition}[section]

\theoremstyle{remark}

\numberwithin{equation}{section}

\newcommand{\vp}{\varphi}
\newcommand{\vep}{\varepsilon}

\newcommand{\Cc}{\mathcal{CC}^{\mathrm{onf}}}

\newcommand{\Cs}{\mathrm{C}^{\mathrm{onf}}}
\newcommand{\R}{\mathbb{R}}

\newcommand{\Pol}{\mathcal{P}}
\newcommand{\Q}{\mathcal{Q}}

\newcommand{\at}{\operatorname{arg}}

\title{On the Configuration Space of Planar Closed Kinematic Chains}
\author{Gerhard Zangerl}
\date{}

\begin{document}

\maketitle

\begin{abstract}
\noindent A planar kinematic chain consists of $n$ links connected by joints. In this work we investigate the space of configurations, described in terms of joint angles, that guarantee that the kinematic chain is closed. We give explicit formulas expressing the joint angles that guarantee closedness by a new set of parameters, the diagonal lengths  (the distances of the joints to the origin) of the closed kinematic chain. Moreover, it turns out that these diagonals are contained in a domain that possesses a simple structure. We expect that the new insight can be applied for several issues such as motion planning for closed kinematic chains or singularity analysis of their configuration spaces. In order to demonstrate practicality of the new method we present numerical examples.   
\end{abstract}

\begin{center}

\end{center}


\section{Introduction}

In this work we investigate  the configuration space of closed planar kinematic chain (CKC) with $n$ links connected by revolute joints in terms of its joint angles. In many fields like robotics computational biology or protein kinematics it is of immense interest to understand the configuration space of a CKC. For instance, in robotics the problem to connect a start, $\alpha_s$, and goal configuration, $\alpha_g$ naturally appears and thus requires knowledge of the configuration space, which is typically a manifold or variety in the ambient space formed by the robots joint variables. The configuration space is even more complicated if additional constraints like obstacle, link-link avoidance, or limited joint angles are included. Two main strategies, probabilistic and geometric approaches, to investigate configuration spaces have been developed so far.\\

Probabilistic methods have been successfully applied for constrained motion planning. They are especially important in practical situations with high dimensions that include complex constraints such as obstacle avoiding. Typically these methods are based on the generation of random configurations in ambient joint space followed by a check up if they approximately satisfy the desired constrains. Repeating this procedure results in a discrete version of the configuration space that is very useful in applications. Probabilistic  methods have been applied in different situations, which can be found in 
\cite{suthakorn2001new,CorSim04,jaillet2013path,SLaV99,SLaV06,zhang2013unbiased,yakey2001randomized},\\

Besides the approaches using randomness other works focused on questions about the geometry and topology of the  configuration spaces of kinematic chains. Insight about the global geometry of configuration spaces is very important in applications. Early discoveries have been made by \cite{KapMil95, HAUKNU98}. In their fundamental work Kapovitch and Milgram established important results about the geometry, which led to novel path planning algorithms. For instance in \cite{JamesRTrinkleJ02a,JamesRTrinkleJ02b} it is used that the configuration space of a CKC consists of two connected components when it possess three long links. An application of this result is that path planning can be done easily for this special kind of CKC's. Also for the more difficult case, when CKC's do not have three long links algorithms were derived in \cite{JamesRTrinkleJ02a,JamesRTrinkleJ02b}. They also developed path planners in the case of $p$ point obstacles in the plane \cite{liu2005toward}. Another geometric approach was recently recognized by  Han, Rudolph and Blumenthal. They discovered that it is very beneficial to describe the configuration space of CKC by different parameters than the joint angles, see \cite{HanRudBlumVal, han2008conv, han2006inverse}. Their idea is to use the length of diagonals from the positions of revolute joints to the origin $O$ as depicted on the right side of Figure \ref{fig:Atan2}. It turns out that for a CKC the length of these diagonals can be computed as solution of a system of linear inequalities, which means that all feasible diagonal lengths can be described by a convex polyhedron that can be handled by methods of linear programming \cite{luenberger2015linear}. Given feasible diagonal lengths, several configurations of the CKC can be constructed, since each link of the chain can be flipped over a diagonal. Thus in \cite{HanRudBlumVal,han2008conv} any configuration can be obtained from a set of diagonals and a vector that represents the choices of flipping, which shows that the configuration space is formed by several copies of the polyhedron given by the system of inequalities. This practically convex structure is very useful for motion planning. In \cite{han2006inverse,han2008conv} paths between CKC with 1000 links are computed very efficiently.\\

\textbf{Contribution of this work:} We develop a new method that explicitly computes configurations of a CKC with $n$ links, which are described by its joint angles. Compared to other methods it does not require linear programming to solve a system of linear inequalities like in \cite{HanRudBlumVal,han2008conv} nor does it rely on probabilistic principles. The developed  method can be used to easily sample configuration space of a CKC and thus is expected to be useful in practical applications.\\

\textbf{Outline of this text:} In section \ref{sec:sampling} we give a mathematical description of a CKC and its configuration space in terms of the diagonals of a CKC. Then the basic algorithm that explicitly describes how configurations of a CKC can be computed from its diagonals is developed in section \ref{sec:sampling}. In section \ref{sec:SDpar} we describe the set of new parameters and show how they can be used to compute a vector of joint angles of  a CKC. Finally,  we give  numerical examples that show the validity of the developed method.

\section{Configuration space}\label{sec:sampling}

To describe the configuration space of a CKC with link lengths $a_1, \dots, a_n$, given as the entries of its vector of link lengths $a^n= (a_1, \dots, a_n)$,  we introduce Cartesian coordinates in two dimensional Euclidean space. Moreover we place one of the links of the CKC so that it is supported by the positive $x$-axis and so that one of its ends coincides with the origin. Without loss of generality we can assume that the link $a_n$ of the chain is fixed in the described manner, see Figure \ref{fig:ClosedChain}. In the following, we identify an angle $\alpha$ with its corresponding point on $S^1$. Furthermore, for $1 \leq k \leq n$ and a vector of angles $\alpha^k := (\alpha_1, \dots, \alpha_{k})\in \left(S^1\right)^k $ we denote by 
\begin{align}
f_{a^n ,k} \colon  \left( S^1\right)^k \rightarrow \R^2, \quad  f_{a^n ,k} (\alpha^k) = \sum_{j=1}^{k} a_j \begin{pmatrix} \cos(\alpha_j) \\
\sin(\alpha_j)
\end{pmatrix}.
\end{align}
the $k$-th endpoint map of a kinematic chain (KC). Additionally, the domain $W_k :=  f_{a^n ,k}\left( (S^1)^k \right)$, which is a circular annulus, will be referred to as the workspace of $k$-th endpoint map. 
We will call $\alpha^{n-1}$ a configuration of the CKC with link lengths $a_1, \dots,a_n$ if it satisfies the closure condition, which means that it is contained in the set
\begin{align}\label{clos_eq}
&& \Cs ({a^n}) = \left\{ \alpha^{n-1} \in \left(S^1\right)^{n-1}  \colon  f_{a^n,n-1}\left( \alpha^{n-1} \right) =
\begin{pmatrix} a_n \\ 0
\end{pmatrix} \right\} = f_{a^n , n-1}^{-1} \left(a_n, 0\right).
\end{align}
If no restrictions on the endpoint map are imposed $\alpha^{n-1}$ will just be called a configuration of the KC with $n-1$ links.  
\begin{figure}
\begin{center}
\begin{tikzpicture}[scale = 0.5]
\draw[thin,->,  -latex] (-1,0)--(18,0);
\draw[thin, ->, -latex] (0,-1)--(0,3);
\draw[thin] (0,0)--( 6,4) -- ( 8,2)--(15, 6) --(16, 0 )  ;
\draw ( 4,4) -- ( 8,4)  ;
\draw ( 6,2) -- ( 11,2)  ;
\draw ( 13,6) -- ( 17,6);
\draw[thin, fill = red] (0,0) circle (4pt);
\draw[thin, fill = red] (6,4) circle (4pt);
\draw[thin, fill = red] (8,2) circle (4pt);
\draw[thin, fill = red] (15,6) circle (4pt);	
\draw[thin, fill = red] (16,0) circle (4pt);
\draw[->, -latex] (1.5,0) arc (0:50:1cm);
\draw[->, -latex] (7,4) arc (0:315:1cm);
\draw[->, -latex] (10,2) arc (0:30:2cm);
\draw[->, -latex] (16,6) arc (0:280:1cm);
\node at (0.9,0.25) {\fontsize{10}{10} $\alpha_1$};
\node at (6,4.5) {\fontsize{10}{10} $\alpha_2$};
\node at (9.3,2.3) {\fontsize{10}{10} $\alpha_3$};
\node at (15,6.5) {\fontsize{10}{10} $\alpha_4$};
\node at (3,2.5) { \rotatebox{35}{$a_1$} };
\node at (7.5,3.1) { \rotatebox{-40}{$a_2$} };
\node at (11.5,4.4) { \rotatebox{40}{$a_3$} };
\node at (15.9,3.1) { \rotatebox{-75}{$a_4$} };
\node at (8,0.4) { \rotatebox{0}{$a_5$} };
\node at (0.8, 2.5)  { \rotatebox{0}{ \fontsize{15}{15}$y$} };
\node at (17.5,-0.5)  { \rotatebox{0}{\fontsize{15}{15}$x$} };
\end{tikzpicture}\caption{A CKC with $n=5$ five links. The link $a_5$ is supported on the positive $x$-axis and one of its ends coincides with the origin}\label{fig:ClosedChain} 
\end{center}
\end{figure}
Furthermore,  the analysis in this work  uses the simple observation that it is sufficient to understand the space \begin{align}\label{eq:circular_configurations}
\Cc({a^n})  =  \left\{   \beta^{n-1} \in \left( S^1\right)^{n-1}  \colon   \| f_{a^n,n-1}\left(\beta^{n-1}	\right) \|^2_2  = a_n^2 \right\},   
\end{align}
in order to describe $\Cs(a^n)$, where $\| \cdot  \|_2 $ denotes the Euclidean norm. From the definition of $\Cc({a^n}) $ it is clear that any configuration $\beta^{n-1} \in \Cc({a^n})$ satisfies that its endpoint 
\begin{alignat*}{2}
f_{a^n,n-1}\left( \beta^{n-1}\right) \in K_{a_n}
\end{alignat*}
lies on the circle $K_{a_n}$ that is centred on the origin and has radius $a_n$. We will say that $\beta^{n-1}$ is closed up to a rotation and call it a \textit{circular configuration} of a CKC. 
Clearly, any circular configuration $\beta^{n-1}$ can be rotated by an angle $\lambda$,   
\begin{alignat*}{2}
\beta^{n-1} + \lambda := \left( \beta_1 + \lambda, \dots , \beta_{n-1} +\lambda  \right),  
\end{alignat*} 
so that $\beta^{n-1} + \lambda  \in C_a $. Thus, if we are able to give an efficient method  to compute the set of solutions to the implicit equation   
\begin{alignat}{2}\label{eq:closed_norm}
\|f_{a^n, n-1} \left(\beta^{n-1}	\right) \|^2_2  = a_n^2, 
\end{alignat}
we also obtain configurations in $\Cs(a^n)$ by the following two step algorithm:
\begin{itemize}
\item[(i)]   Compute a circular configuration  $\beta^{n-1} \in \Cc({a^n})$
\item[(ii)]  Determine $\lambda$ such that  $\alpha^{n-1} =  \beta^{n-1} + \lambda  \in \Cs({a^n})$ 
\end{itemize}
Once a circular configuration is obtained step (ii) is a rather simple task. Therefore, in the following we will focus on the solution of step (i). This step is based on the fact that the trigonometric equation \eqref{eq:closed_norm}, which in its expanded form is given as    
\begin{align}\label{eq:CKCsquared_expanded}
\sum_{i=1}^{n-1} a_i^2 + 2 \sum_{i<j}^{n-1} a_i a_j 
\cos(\beta_i -\beta_j) = a_n^2 , 
\end{align}
allows for some kind of backwards substitution, see section \ref{sec:CC}. By the preimage theorem we know that the set of all circular configurations of a CKC with $n$ links satisfying  \eqref{eq:CKCsquared_expanded} is a manifold of dimension $n-2$, whenever $a_n^2$ is a regular value of the map $g \left(\beta^{n-1} \right)  := \| f_{a, n-1}\left( \beta^{n-1} \right)    \|^2_2$. In all other cases the space $\Cc_{a}$ may have singular points.

\subsection{Mathematical tools and notations}\label{sec:noation}

Surprisingly, the trigonometric equation \eqref{eq:CKCsquared_expanded} can be rearranged into an equation of the same type but with one joint angle less appearing on its left hand side. For the computations we use  that a linear combination of sine and cosine functions can be written as      
\begin{align}\label{eq:addsincos}
a \sin\left( x \right) + b \cos\left( x\right) = c \sin\left(x + \varphi(a,b)  \right), 
\end{align}
where $c = \sqrt{a^2 + b^2}$  \text{and}  $\varphi\left(a,b\right)  = \at \left(a, b\right) $
is the function described in Figure \ref{fig:Atan2}.
\begin{figure}
\centering 
\begin{tikzpicture}[scale = 0.7]
\begin{axis}[ 
axis x line  = middle, 
axis y line = middle, xmin = -2, xmax = 2, ymin = -0, ymax = 2, 
xtick = \empty,  
ytick = \empty,
grid = major, enlargelimits = 0.2]
\draw[->, -latex, line width = 1pt] (axis cs: 0,0)--(axis cs:-2,2);
\draw[->, -latex] (axis cs:1.4,0) arc (0:114:2cm);
\node at (axis cs:-1.8, 2.2) {\fontsize{15}{15} $P(a|b)$ };
\node at (axis cs: 0.35,0.5)  { \colorbox{white}{\fontsize{13}{13} $\vp(a,b)$}};
\node at (axis cs:2.4, -0.2) {\fontsize{15}{15} $x$};
\node at (axis cs:0.3, 2) {\fontsize{15}{15} $y$};
\draw[dashed] (axis cs: -2,0)--(axis cs: -2,2);
\draw[dashed] (axis cs:0,2)--(axis cs: -2,2);
\draw[thin, fill = red] (axis cs: 0,0) circle (3pt);
\end{axis}
\end{tikzpicture}\hspace{1cm}
\begin{tikzpicture}[scale = 0.7]
\begin{axis}[ 
axis x line  = middle, 
axis y line = middle, xmin = 0, xmax = 5, ymin = -0, ymax = 3, 
xtick = \empty,  
ytick = \empty,
grid = major, enlargelimits = 0.2]
\draw (axis cs: 0, 0)--(axis cs: 1, 3);
\draw (axis cs: 2.5, 3)-- (axis cs:4, 2);
\draw (axis cs: 4, 2)--(axis cs: 5, 0); 
\draw (axis cs: 0, 0)--(axis cs: 4, 2);
\draw (axis cs: 0, 0)--(axis cs: 2.5, 3);
\draw[dashed] (axis cs: 1, 3)--(axis cs: 2.5, 3);
\draw[fill = red] (axis cs: 0 ,0 )   circle (3pt);
\draw[fill = red] (axis cs: 1 ,3 )   circle (3pt);
\draw[fill = red] (axis cs: 4 ,2 )   circle (3pt);
\draw[fill = red] (axis cs: 5 ,0 )   circle (3pt);
\draw[fill = red] (axis cs: 2.5 ,3) circle (3pt);
\node at (axis cs:2.5, -0.4 ) {\fontsize{15}{15}$L_{n-1} =  a_n$};
\node at (axis cs:3.3, 2.7 ) {\fontsize{15}{15} \rotatebox{-40}{$a_{n-2}$}};
\node at (axis cs:4.8, 1 ) {\fontsize{15}{15} \rotatebox{-70}{$a_{n-1}$}};
\node at (axis cs:0.4, 2.3 ) {\fontsize{15}{15} \rotatebox{80}{$L_1 =  a_{1}$}};
\node at (axis cs:2.2, 1.5 ) {\fontsize{15}{15} \rotatebox{40} { \fontsize{15}{15} $L_{n-2}  $}};
\end{axis}
\end{tikzpicture}\caption{Left: The function $\at$ gives the angle between the $x$-axis and the vector from the origin to $P(a|b)$. Right: A circular configuration with endpoint $(a_n, 0)$. The picture shows anchored diagonals of the CKC}\label{fig:Atan2}
\end{figure}
In order to  achieve a compact presentation of the results that will follow it is important to introduce abbreviations. For this purpose consider   
\begin{align*}
2a_{n-1} \sum_{j=1}^{n-2} a_j \cos\left(\beta_{n-1} - \beta_j  \right) + 2 \sum_{i<j}^{n-2} a_i a_j  \cos(\beta_i -\beta_j) + \sum_{i=1}^{n-1} a_i^2 =  a_n^2     \end{align*}
which is an equivalent form of \eqref{eq:CKCsquared_expanded} that is obtained by fixing an index to be $n-1$ and rearranging remaining terms. Finally, using trigonometric summation formulas we arrive at    
\begin{align}\label{eq:squaredrearranged}
2 a_{n-1} \sin\left( \beta_{n-1} \right) \sum_{j=1}^{n-2} a_j \sin\left( \beta_j \right)  &+ 2a_{n-1} \cos\left( \beta_{n-1} \right) \sum_{j=1}^{n-2} a_j \cos\left( \beta_j \right) +  \sum_{i=1}^{n-1} a_i^2   \nonumber \\ 
&  +2\sum_{i<j}^{n-2} a_i a_j \cos\left( \beta_i - \beta_j \right) =  a_n^2.  
\end{align}
In the last expression addition formula \eqref{eq:addsincos} can be applied, which motivates the following abbreviations:\\

For a CKC with link lengths vector $a^n = (a_1, a_2, \dots, a_n)$ for $1\leq k \leq n-1$ we set $\beta^{k} := (\beta_1, \dots, \beta_{k})$ and we abbreviate the $x$- and $y$-coordinates of the $k$-th endpoint map $f_{a^n,k}$ by    
\begin{alignat}{2}\label{eq:coordinates}
X_k(\beta^k):=  \sum_{i=1}^k a_i \cos(\beta_i)~~ \text{and }~~ Y_k(\beta^k) :=  \sum_{i=1}^k a_i \sin(\beta_i).
\end{alignat}
Using this notation we define 
\begin{alignat*}{2}
L_{k}\left(\beta^{k} \right) 	&:= \sqrt{X_k(\beta^k)^2  + Y_k(\beta^k)^2  }  =  \sqrt{ \sum_{i=1}^{k} a_i^2 +     2 \sum_{i< j}^{k} a_i a_j \cos\left(\beta_i - \beta_j\right)},
\end{alignat*}
which we refer to as the $k$-th diagonal length that  by definition can also be written as  $ L_k(  \beta^{k} ) =  || f_{a^n, k}(\beta^k)||$. Moreover, consider the angle  \begin{alignat*}{2}
\Phi_k \left(\beta^{k}\right)	&:= \varphi \left(  Y_k\left( \beta^k \right),  X_k\left( \beta^k \right)     \right), \end{alignat*}
which naturally appears, when applying formula \eqref{eq:addsincos} for equation \eqref{eq:squaredrearranged}. We set $\beta^1 := \beta_1$. We point out here that $L_1\left( \beta^1 \right) = a_1$, $L_{n-1}\left( \beta^{n-1} \right) = a_n$. Additionally, we denote by $L^{n-3} = (L_2, \dots L_{n-2}) \in \R^{n-3}$ the vector of (variable) diagonals for a KC. Moreover, we use the notation 
\begin{align}\label{eq:RmincmaxA}
\mathrm{R}_{k}^{\min}   := 0 \vee  \max_{1\leq i\leq k}\left(2a_i-\sum_{j=1}^ka_j\right) ~ \text{and}~ \mathrm{    R}_{k}^{\max} :=      \sum_{j=1}^ka_j  ~~\text{for}~1\leq i\leq k, 
\end{align} 
where $0 \vee a := \max\left\{ 0, a \right\}$ for $a\in \R$. These quantities denote the inner and outer radius of the workspace $W_k$ of the KC with links $a^k = (a_1, \dots, a_k)$. Note that in computations carried out in section \ref{sec:CC} arguments of the defined quantities will frequently be omitted for simplicity and that with the notations equation \eqref{eq:squaredrearranged} can be rewritten as
\begin{align}\label{eq:SQshort}
2 a_{n-1} \sin\left( \beta_{n-1} \right) Y_{n-2}\left( \beta^{n-2} \right)  &+ 2a_{n-1} \cos\left( \beta_{n-1} \right) X_{n-2}\left( \beta^{n-2} \right) +  L_{n-2}^2 =  L_{n-1}^2 - a_{n-1}^2.  
\end{align}

\subsection{Diagonal Space and Circular Configurations}\label{sec:CC} 

We give a new method to obtain solutions to equation \eqref{eq:CKCsquared_expanded}. In the proof of Theorem \ref{thm:sampling} a procedure how to obtain such solutions by reducing the length of the CKC step by step is described. For $k \geq 1$ the solution method involves the choice of a real value in a domain that is given by an inequality that involves the diagonals of the CKC. Before we state the main theorem we give the following definition.      

\begin{definition}[Diagonal Space]\label{def:DiagSpace} For a KC with a vector of links  $a^n= ( a_1,\dots, a_n) $ and recall that $L^{n-3} \in \mathbb{R}^{n-3}_{\geq 0}$. We denote the diagonal space of the KC as the set \begin{align}
&\mathcal{DS}({a^n}):= \left\{  L^{n-3}  \colon L_{n-k-1} \in   
\Bigl[\left|L_{n-k}-a_{n-k}\right|,L_{n-k}+a_{n-k}\Bigr] \cap  
\Bigl[\mathrm{R}_{n-k-1}^{\min}, \mathrm{    R}_{n-k-1}^{\max}\Bigr], ~1\leq k \leq n-3 \right\}.   
\end{align}
\end{definition}

\begin{theorem}[Computation of circular configurations]\label{thm:sampling} Every circular configuration
\(\beta^{n-1}\in  \Cc(a^n)\)
has a diagonal vector contained in
\(\mathcal{DS}(a^n)\).
Conversely, for every
\(L^{n-3}\in\mathcal{DS}(a^n)\),
there exists a circular configuration having these diagonal lengths. Moreover, the the angle $\beta_k$ is related to the angles $\beta^{k-1}$ by the formula 
\begin{align}\label{eq:diagangle}
2 a_k L_{k-1}  \sin\left( \beta_k + \Phi_{k-1}\left( \beta^{k-1}\right)  \right) =  L_k^2  - a_k^2 - L_{k-1}^2,   
\end{align}
for $1 \leq k \leq n-1$, whenever $\Phi_{k-1}$ is defined. If not, the relation between $\beta_k$ and $\beta^{k-1}$ is given by \eqref{eq:SQshort}, where $n-1$ and $n-2$ are replaced by $k$ and $k-1$. Note that for $k=1$ we set $L_0 := 0$ and $\Phi_0 :=0$ so that $\beta_1$ is an arbitrary angle. 
\end{theorem}

\begin{proof} Assume $\beta^{n-1}$ is a circular configuration and thus solves equation \eqref{eq:SQshort}. We will manipulate equation \eqref{eq:SQshort} to show that the vector of variable diagonal lengths $L^{n-2}$ is indeed an element of $\mathcal{DS}({a^n})$. We start by showing  that the diagonal $L_{n-2}$ is contained  in $[|L_{n-1}-a_{n-1}|, L_{n-1}+a_{n-1}] \cap [ \mathrm{    R}_{n-2}^{\min} , \mathrm{    R}_{n-2}^{\max}]$. For the remaining diagonals the analogous statement for $L_{n-k-1}$ then follows inductively by applying the same arguments.\\

Although it is very tempting to apply addition formula \eqref{eq:addsincos} to this equation we have to  deal with special cases before we can do so.\\

$X_{n-2} \neq 0$ but $Y_{n-2} = 0$. In this case, according to \eqref{eq:SQshort} we obtain the equation  \begin{align*}
& 2 a_{n-1} \cos(\beta_{n-1} )  X_{n-2} + L_{n-2}^2  =  L_{n-1}^2 - a_{n-1}^2.  
\end{align*}
solving for $\beta_{n-1}$ yields  
\begin{align*}
& \cos(\beta_{n-1} )     =  \frac{ L_{n-1}^2 - L_{n-2}^2- a_{n-1}^2 }{ 2 a_{n-1}  X_{n-2}  }.  
\end{align*}
The Latter equation can be solved for $\beta_{n-1}$ iff the square of its right side is lower or equal to one. Since in the considered case $ X_{n-2}^2 = L_{n-2}^2$,  this leads to the inequality   
\begin{align}\label{ieq:diagonals}
\left(   L_{n-1}^2 - L_{n-2}^2- a_{n-1}^2 \right)^2 \leq  4 a_{n-1}^2 L_{n-2}^2, 
\end{align}
that is satisfied for $L_{n-2} \in [ |L_{n-1} -a_{n-1}|,  L_{n-1} +a_{n-1} ]$. Additionally, since $L_{n-2}$ is a diagonal length also $L_{n-2} \in [\mathrm{R}_{n-2}^{\min}, \mathrm{    R}_{n-2}^{\max}]$ has to be satisfied. The case when $X_{n-2} = 0$ and $Y_{n-2} \neq 0$ is analogously.\\       

The case $ L_{n-2} = X_{n-2} = Y_{n-2} = 0$ can only  occur for $a_n = L_{n-1} = a_{n-1}$. Here $\beta_{n-1}$ is an arbitrary value, and $\mathrm{R}_{n-2}^{\min} = 0$ since $L_{n-2} = 0 $  is a diagonal of the circular configuration. Also $0 = L_{n-2} \in [ |L_{n-1} -a_{n-1}|,  L_{n-1} +a_{n-1} ] = [0, 2 a _n ]$ is satisfied.\\ 

If $X_{n-2} \neq 0 $ and  $Y_{n-2} \neq 0 $ we can apply the addition theorem to formula \eqref{eq:SQshort} and obtain   
\begin{align*}
& 2 a_{n-1} \sqrt{  X_{n-2}^2 + Y_{n-2}^2 } \, \sin\left( \beta_{n-1} + \Phi_{n-2}\left(\beta^{n-2} \right)  \right) + L_{n-2}^2  = L_{n-1}^2 - a_{n-1}^2,   
\end{align*}
which can be rewritten as 
\begin{align*}
& 2 a_{n-1} L_{n-2} \sin\left( \beta_{n-1} + \Phi_{n-2}\left(\beta^{n-2} \right)  \right) + L_{n-2}^2  = L_{n-1}^2 - a_{n-1}^2.   
\end{align*}
Repeating the arguments outlined above We conclude  that this equation can only have a solution if \eqref{ieq:diagonals} is satisfied and thus also in this case $L_{n-2}  \in  [ |L_{n-1} -a_{n-1}|,  L_{n-1} +a_{n-1} ]$ and $L_{n-2} \in [\mathrm{R}_{n-2}^{\min}, \mathrm{    R}_{n-2}^{\max}]$. Thus, repeating the arguments it follows that $L^{n-3} \in \mathcal{DS}(a^n)$.\\
       
Conversely, if $L^{n-3} \in \mathcal{DS}(a^n)$ then  equation \eqref{eq:diagangle} can be solved for $\beta_k$ if $\Phi_{k-1}$ is defined. If not, equation \eqref{eq:SQshort} can be used alternatively, with $n-1$, $n-2$ replaced by $k$ and $k-1$. Repeating this procedure gives a circular configuration $\beta^{n-1}$. 
\end{proof}

Theorem  \ref{thm:sampling} outlines a procedure how the entries of elements $L^{n-3} \in \mathcal{DS}(a^n)$ can be determined recursively. This leads to a sampling strategy for the whole diagonal space. Moreover, equation \eqref{eq:diagangle} explicitly describes how angles can be obtained from  the diagonals. We shortly summarize this in more detail in order to emphasize the progress of this work in comparison with the existing literature: 
\begin{itemize}
\item The explicit sampling procedure for $\mathcal{DS}(a^n)$ is an advance over the previous methods  \cite{han2006inverse,han2008conv}, where diagonals are obtained by linear programming. 
\item Looking more closely at the definition of $\mathcal{DS}(a^n)$, one finds that the diagonal space is the intersection of a polytope and a cuboid.  More precisely, if 
\begin{align}\label{eq:Pol}
\Pol(a^n):=\left\{(L_2,\dotsc,L_{n-2})\in \mathbb{R}^{n-3}_{\geq 0}   \colon L_{n-k-1} \in   \Bigl[\left|L_{n-k}-a_{n-k}\right|,L_{n-k}+a_{n-k}\Bigr]\right\},
\end{align}
denotes the polytope defined by nested intervals, then
\begin{align}\label{eq:Cub}
 \mathcal{DS}(a^n) = \Pol(a^n) \cap \Q(a^n), ~~\text{with}~~ \Q(a^n) :=  \prod_{1\leq k \leq n-3} \Bigl[ \mathrm{ R}_{n-k-1}^{\min}, \mathrm{  R}_{n-k-1}^{\max}\Bigr].
\end{align}
This is a new representation of the diagonal space. 
\end{itemize}
Although the derived sampling procedure already provides a useful description of the diagonal space $\mathcal{DS}(a^n)$,  we will investigate $\mathcal{DS}(a^n)$ further in  section \ref{sec:SDpar}.

\section{Circular Configurations from Diagonals}\label{sec:diagspace}

In this section, we take a closer look at equation \eqref{eq:diagangle}, which is solved recursively to recover the joint angles from a feasible diagonal vector. By Theorem \ref{thm:sampling}, a vector of joint angles $\beta^{n-1}$ can be computed from $L^{n-3}\in\mathcal{DS}(a^n)$ using
\begin{align}\label{eq:diagangleB}
2 a_k L_{k-1}  \sin\left( \beta_k + \Phi_{k-1}\left( \beta^{k-1}\right)  \right) =  L_k^2  - a_k^2 - L_{k-1}^2 
\end{align}
for $1\leq k\leq n-1$. Since $L_0=0$ and $L_1=a_1$, equation \eqref{eq:diagangleB} is identically satisfied for $k=1$, which means that the entire configuration can be rotated  around the origin by an arbitrary angle.

For $2\leq k\leq n-1$ we assume that $L_{k-1}>0$, which implies that equation \eqref{eq:diagangleB} can be rearranged 
\begin{alignat}{2}\label{eq:beta_epsilon}
\beta_k^{\vep_k}
=
\pi\vep_k
+
(-1)^{\vep_k}S_k
-
\Phi_{k-1}\left(
\beta^{k-1}
\right),
\qquad
\vep_k\in\{0,1\},
\end{alignat}
where we abbreviate
\begin{align*}
S_k
:=
\arcsin\left(
\frac{L_k^2-a_k^2-L_{k-1}^2}{2a_kL_{k-1}}
\right).
\end{align*}
Here, $\vep_k$ records which preimage is selected at the $k$-th step and all angles are understood modulo $2\pi$. Note that if the argument of the arcsine is equal to $1$ or $-1$, the two expressions in \eqref{eq:beta_epsilon} agree modulo $2\pi$; hence, there is only one distinct preimage at that step.\\

It remains to consider the case $L_{k-1}=0$, which implies that $L_k=a_k$. In this situation and equation \eqref{eq:diagangleB} is again identically satisfied and therefore $\beta_k$ can be chosen arbitrarily, as depicted in Figure \ref{fig:ClosedChain_b}. 

We can now describe the complete preimage of any $L^{n-3}\in\mathcal{DS}(a^n)$ under the map that assigns to each circular configuration its diagonal vector. For every $1\leq k \leq n-1$ that satisfies $L_{k-1}>0$ compute its (possibly coinciding) preimages by \eqref{eq:beta_epsilon}. At every step with $L_{k-1}=0$, the angle $\beta_k$ is chosen arbitrarily. Conversely, every sequence of choices made in this way produces a circular configuration with the prescribed diagonal vector.

For $L_{k-1}>0$, the choice of $\vep_k$ has a simple geometric interpretation: whenever the two preimages are distinct, they correspond to the two placements of the link $a_k$ obtained by reflecting it across the line that is supported by $L_{k-1}$, as illustrated in Figure \ref{fig:ClosedChain_b}. Equivalently, they are the two possible orientations of the corresponding triangle in the terminology of \cite{han2008conv}.
\begin{figure}[h!]
\center 
\begin{tikzpicture}[scale = 0.9]

\draw[thin,->,  -latex] (-2,0)--(4,0);
\draw[thin, ->, -latex] (0,-1)--(0,2);
\draw[thin] (0:0cm)--(170:2cm)--(-1,2);
\draw[thin, dashed] (-1,2)--(1,2.5);
\draw[thin] (1,2.5)--(2.5, 2.5)--(0,0)--(3,0.5)--(3.5,-0.5);
\draw[thin, dashed] (3.5,-0.5)--(2.5,-1);
\draw[thin] (2.5,-1)--(1.5,-0.5)--(1,0);
\draw[fill = red ]    (0,0) circle (2pt);
\draw[fill = red ]    (170:2cm) circle (2pt);
\draw[fill = red ]    (-1,2) circle (2pt);
\draw[fill = red ]    (1,2.5) circle (2pt);
\draw[fill = red ]    (2.5, 2.5) circle (2pt);
\draw[fill = red ]    (3,0.5) circle (2pt);
\draw[fill = red ]    (3.5,-0.5) circle (2pt);
\draw[fill = red ]    (2.5,-1) circle (2pt);
\draw[fill = red ]    (1.5,-0.5) circle (2pt);
\draw[fill = red ]    (1,0) circle (2pt);
\node at ( -1,0.35) { \fontsize{10}{10} \rotatebox{-5}{$a_{1}$} };
\node at ( -1.6,1.5) { \fontsize{10}{10} \rotatebox{65}{$a_{2}$} };
\node at ( 1.8,2.7) { \fontsize{10}{10} \rotatebox{0}{$a_{k-2}$} };
\node at ( 1,1.4) { \fontsize{10}{10} \rotatebox{45}{$a_{k-1}$} };
\node at ( 1.6,0.5) { \fontsize{10}{10} \rotatebox{10}{$L_{k} = a_{k}$} };
\node at ( 1.8,-0.3) { \fontsize{10}{10} \rotatebox{0}{$a_{n-1}$} };
\node at ( 1.4,-0.8) { \fontsize{10}{10} \rotatebox{0}{$a_{n-2}$} };
\node at ( 0.5,-0.2) { \fontsize{10}{10} \rotatebox{0}{$a_{n}$} };
\end{tikzpicture}\hspace{1cm}
\begin{tikzpicture}[scale = 0.9]
\coordinate (A) at (50:5cm);
\coordinate (B) at (12:5cm);
\coordinate (R) at (30:3cm);
\coordinate (Q) at (4,1.6);
\draw[thin,->,  -latex] (-2,0)--(1,0);
\draw[thin, ->, -latex] (0,-1)--(0,2);
\draw[thin] (0:0cm)--(150:2cm);
\draw[thin] (150:2cm)--(110:3cm);
\draw[thin, dashed] (110:3cm)--(0.5,2.5);
\draw[thin] (0.5,2.5)--(30:3cm);
\draw[thin, dashed] (0,0)--(30:6.5cm);
\draw[thin] (30:3cm) --(50:5cm);
\draw[thin] (30:3cm) --(10:5cm);
\draw[thin, dashed] (1,1.5) --(5,1.5);
\draw[thin, dashed]   (10:5cm)--(50:5cm);
\pic [draw, thin, ->,-latex, angle eccentricity=1.5, angle radius = 0.6cm] {angle = Q--R--A};
\pic [draw, thin, ->,-latex, angle eccentricity=1.5, angle radius = 0.4cm] {angle = Q--R--B};
\draw[fill = red ]    (150:2cm) circle (2pt);
\draw[fill = red]     (0:0cm)   circle (2pt);
\draw[fill = red]     (30:3cm)  circle (2pt);
\draw[fill = red]     (50:5cm)  circle (2pt);
\draw[fill = red]     (10:5cm)  circle (2pt);
\draw[fill = red]     (0.5,2.5) circle (2pt);
\draw[fill = red]     (110:3cm) circle (2pt);
\node at ( -1.3,0.5) { \fontsize{10}{10} \rotatebox{5}{$a_{1}$} };
\node at ( -1.3,1.5) { \fontsize{10}{10} \rotatebox{5}{$a_{2}$} };
\node at ( 1.5,2.2) { \fontsize{10}{10} \rotatebox{-25}{$a_{k-1}$} };
\node at ( 2.6,2.5) { \fontsize{10}{10} \rotatebox{75}{$a_{k}$} };
\node at ( 1,1) { \fontsize{10}{10} \rotatebox{0}{$l$} };
\end{tikzpicture}
\caption{Left: Situation where $L_{k-1} = 0$ and thus $L_k=a_k$, which means that $\beta_k$ can be chosen arbitrarily. Right: Geometric interpretation of the reconstruction from diagonal lengths. For $L_{k-1}>0$, the two distinct preimages in \eqref{eq:beta_epsilon} are related by a reflection of the link $a_k$ across the diagonal $l$. If $L_{k-1}=0$, the corresponding triangle degenerates and the angle $\beta_k$ can be chosen arbitrarily.}
\label{fig:ClosedChain_b}
\end{figure}

\section{Numerical simulations}\label{sec:Appcon}

In this section, we provide numerical examples illustrating the methods developed in this work. We consider CKCs with five and six links. For CKCs with five links, we depict the spaces $\mathcal{DS}(a^n)$. We randomly select elements of the diagonal spaces, compute corresponding circular configurations using Theorem \ref{thm:sampling}, and depict them in Figures \ref{fig:Configs5R} and \ref{fig:configs6R}.

In the numerical reconstruction, equation \eqref{eq:beta_epsilon} is evaluated only when $L_{k-1}>0$. At such a step, either of the two preimages may be selected. If $L_{k-1}=0$, the quotient in the definition of $S_k$ is not evaluated; instead, $\beta_k$ is assigned an arbitrary value. In our computations, we choose $\beta_k=0$ in this case. This convention selects one representative from the continuous family of preimages and does not affect the prescribed diagonal lengths. In floating-point computations, the condition $L_{k-1}=0$ is tested using a suitable numerical tolerance.

\subsection{CKCs with five and six links}

\begin{figure}[htb!]
\centering
\includegraphics[width=0.4\textwidth]{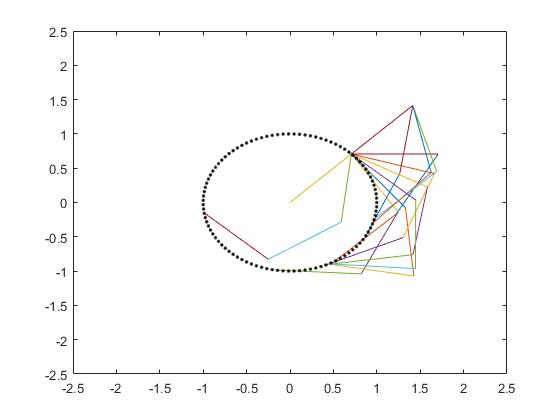}\hspace{1cm}
\includegraphics[width=0.4\textwidth]{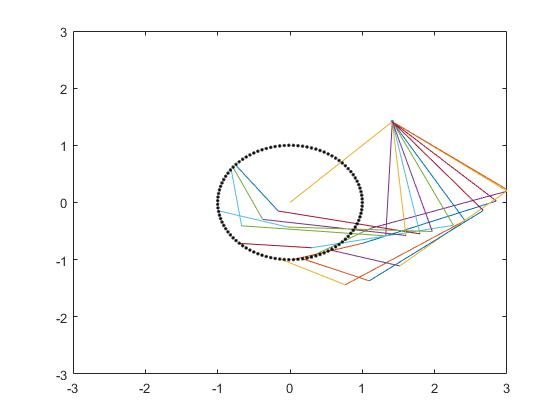}
\caption{Circular configurations of two CKCs with five links. Left: Ten random circular configurations are depicted for the CKC with link lengths equal to one. Right: Ten random circular configurations are depicted for the CKC with link lengths $2,2,2,1,1$.  }\label{fig:Configs5R}
\end{figure}

\begin{figure}[htb!]
\centering
\includegraphics[width=0.4\textwidth]{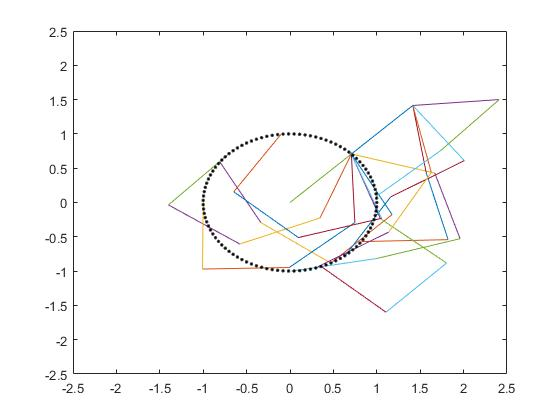}\hspace{1cm}
\includegraphics[width=0.4\textwidth]{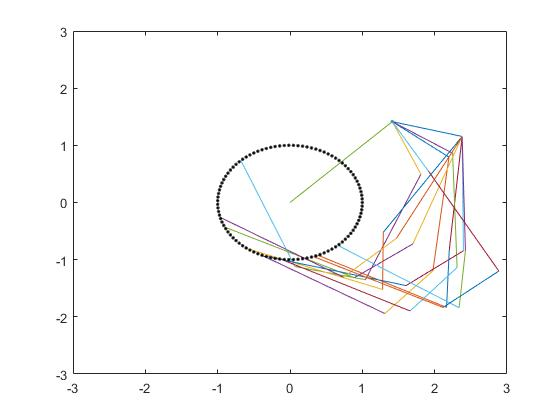}
\caption{Circular configurations of two CKCs with six links. Left: Ten random circular configurations are depicted for the CKC with link lengths equal to one. Right: Ten random circular configurations are depicted for the CKC with link lengths $2,1,2,1,2,1$. }\label{fig:configs6R}
\end{figure}

First, we consider CKCs with $n=5$ links. The domains $\mathcal{DS}(a^5)$ are depicted for two CKCs in Figure \ref{fig:DSspace}. Figure \ref{fig:Configs5R} shows $10$ random circular configurations for these CKCs, where the length of the last link is equal to the radius of the depicted circle. At every step with $L_{k-1}>0$, the first expression in \eqref{eq:beta_epsilon} was selected by setting $\vep_k=0$. At a step with $L_{k-1}=0$, we used the convention $\beta_k=0$ described above.
\begin{figure}[htb]
\centering
\begin{tikzpicture}[scale = 0.7]
\draw[line width = 0.5pt, gray] (0,0) grid (6,6);
\draw[->, -latex, line width = 0.5pt] (0,0) -- (7,0);
\draw[->, -latex, line width = 0.5pt] (0,0) -- (0,7);
\draw[line width = 1pt, fill = red, opacity = 0. 5] (0,2) -- (2, 0) -- (4, 2)--(4, 6 )--(0,2);
\draw[line width = 1pt, fill = red, opacity = 0. 5] (0,0) -- (6, 0) -- (6, 4)--(0, 4 )--(0,0);
\draw[line width = 1.5pt] (0,2) -- (2, 4) -- (4, 4)--(4, 2 )--(2,0)--(0,2);
\node at (6.5,0.5 ) {$L_3$}; 
\node at (0.5, 6.5) {$L_2$}; 
\node at (3.5, 4.6) {  \colorbox{red!50!white}{$\mathcal{P}$} }; 
\node at (5, 3.5) {  \colorbox{red!50!white}{$\mathcal{Q}$} }; 
\node at (2, -0.5) {  \colorbox{white}{$1$} }; 
\node at (-0.5, 2) {  \colorbox{white}{$1$} }; 
\node at (2, 2) {  \colorbox{red!75!white}{$\mathcal{P} \cap \mathcal{Q}$ } }; 
\end{tikzpicture}
\hspace{1.5cm}
\begin{tikzpicture}[scale = 0.5]
\draw[ line width = 0.5pt, gray] (0,0) grid (12,8);
\draw[->, -latex, line width = 0.5pt] (0,0) -- (13,0);
\draw[->, -latex, line width = 0.5pt] (0,0) --(0,10);
\draw[line width = 1pt, fill = red, opacity = 0. 5] (4,8) -- (0,4) -- (4, 0) ;
\draw[line width = 1.5pt] (4,8) -- (0,4) -- (4, 0)--(4,8) ;
\draw[line width = 1pt, fill = red, opacity = 0. 5] (0,0) -- (12, 0) -- (12, 8)--(0,8 )--(0,0);
\draw[line width = 1.5pt] (4,8) -- (0,4) -- (4, 0)--(4,8) ;
\node at (12.5,0.5 ) {$L_3$}; 
\node at (0.5, 9.5) {$L_2$}; 
\node at (2.5, 3) {  \colorbox{red!75!white}{$=\mathcal{P}$} }; 
\node at (7, 3.5) {  \colorbox{red!50!white}{$\mathcal{Q}$} }; 
\node at (4, -0.5) {  \colorbox{white}{$2$} }; 
\node at (-0.5, 4) {  \colorbox{white}{$2$} }; 
\node at (2.5, 4) {  \colorbox{red!75!white}{$\mathcal{P} \cap \mathcal{Q}$ } }; 
\end{tikzpicture}
\caption{Left: Right: The Domain $\mathcal{DS}(a^5)$ for a CKC with five links of length one as intersection of a trapezoidal domain  with a rectangle. Right: The diagonal space of the CKC with links $(2,2,2,1,1)$. Note that on the right side $ \mathcal{DS}(a^5) =   \mathcal{P}(a^5) \cap \mathcal{Q}(a^5) = \mathcal{P}(a^5)$.    }\label{fig:DSspace}
\end{figure}
Finally, we give examples for CKCs with six links. Figure \ref{fig:configs6R} shows random configurations obtained by selecting the first expression in \eqref{eq:beta_epsilon} at every step with $L_{k-1}>0$, that is, by setting $\vep_k=0$. At a step with $L_{k-1}=0$, we again set the freely selectable angle $\beta_k$ equal to zero.

\section{Conclusion and future work}

We have developed a new method for computing configurations in terms of joint angles of a CKC by a systematic procedure. Our approach does not require the solution of a system of linear inequalities by linear programming, nor does it rely on probabilistic methods. Numerical examples demonstrate the validity of the proposed work. We expect that the described method can be useful in tasks such as motion planning for CKCs. We anticipate that it will be an interesting approach for future work to further investigate the diagonal space $\mathcal{DS}(a^n)$. Here it would be interesting to investigate how special designs for CKCs are reflected in the geometry of the diagonal space.




\end{document}